\documentclass[english,oneside]{amsart}
\usepackage[T1]{fontenc}
\usepackage[latin1]{inputenc}
\usepackage{amssymb}
\makeatletter

\newcommand{\cl}[1]{\mathcal{#1}}
\newcommand{\bb}[1]{\mathbb{#1}}
 \theoremstyle{plain}
 \newtheorem{thm}{Theorem}[section]
 \numberwithin{equation}{section} 
 \numberwithin{figure}{section} 
 \theoremstyle{plain}
 \newtheorem{lem}[thm]{Lemma} 
 \theoremstyle{plain}
 \newtheorem{prop}[thm]{Proposition} 
 \theoremstyle{plain}
 \newtheorem{cor}[thm]{Corollary} 
 \theoremstyle{definition}
  \newtheorem{example}[thm]{Example}
\theoremstyle{remark}
\newtheorem{rem}[thm]{Remark}

\makeatother
\begin{document}

\title{On the Classification of Multi-Isometries}

\author{H. Bercovici, R.G. Douglas, and C. Foias}

\thanks{The authors were supported in part by grants from the
National Science Foundation}
\begin{abstract}
We consider the classification, up to unitary equivalence, of commuting
$n$-tuples $(V_{1},V_{2},\dots,V_{n})$ of isometries on a Hilbert
space. As in earlier work by Berger, Coburn, and Lebow, we start by
analyzing the Wold decomposition of $V=V_{1}V_{2}\cdots V_{n}$, but
unlike their work, we pay special attention to the case when $\ker V^{*}$
is of finite dimension. We give a complete classification of $n$-tuples
for which $V$ is a pure isometry of multiplicity $n$. It is hoped
that deeper analysis will provide a classification whenever $V$ has
finite multiplicity. Further, we identify a pivotal operator in the
case $n=2$ which captures many of the properties of a bi-isometry.
\end{abstract}
\maketitle \emph{The present work was initiated in 1999, during the
memorial conference in honor of B\'ela Sz.-Nagy. This paper is a
token of his perennial presence in mathematics.}

\section{Introduction}

Much work in operator theory, particularly the model theory
\cite{key-2,Sz-CBMS} of B. Sz.-Nagy and the third author, relies on
a good understanding and classification of isometric operators on
Hilbert space. This understanding was initiated by J. von Neumann in
a foundational paper on operator theory \cite{key-4} where he
demonstrates the decomposition of an isometry into a direct sum of a
unitary and a unilateral shift. This decomposition was later
rediscovered by H. Wold \cite{key-5}, who made it the cornerstone of
prediction theory for stationary random processes. The deep
relationship between harmonic analysis and shifts on Hilbert space
was discovered then by A. Beurling \cite{key-7}. The explosive
development of operator theory and harmonic analysis which followed
from these discoveries and Sz.-Nagy's dilation theory \cite{Sz1,Sz2}
continues to this day.

By contrast, the understanding and classification of commuting pairs
or, more generally, commuting families of isometries is very partial.
A set of unitary invariants for finite sets of commuting isometries
was found by C. Berger, L. Coburn and A. Lebow \cite{key-1}; these
invariants will also be considered in Section 2 below in somewhat
more detail. Berger, Coburn and Lebow use these invariants in demonstrating
that there are infinitely many nonisomorphic $C^{*}$-algebras generated
by pairs of commuting isometries. Other authors
\cite{ag-cl-douglas,douglas,do-young,gh-mand,salas}
also demonstrate the great variety of families of commuting isometries,
thus showing how difficult the classification problem can be.

In this note we review the unitary invariants of $n$-tuples of commuting
isometries, originally introduced in \cite{key-1} (cf. Theorems 2.1
and 2.8), we produce a smaller set of invariants (Theorem 2.3), and
we determine some of the necessary conditions these invariants must
satisfy (Proposition 2.4 and Theorem 2.5). The most specific results
in this section pertain to the case $n=3$. In this case, the possibility
of finding triples of commuting isometries is related with the existence
of invariant subspaces for a specific contraction (cf. Theorem 2.5
and Corollary 2.6). It is difficult to conjecture the natural extension
of these results to $n\ge 4$. In Section 3 we obtain a deeper understanding of bi-isometries relating their structure to the model theory for contractions. Finally, in Section 4, we provide a complete
classification for a certain class of irreducible $n$-tuples of commuting
isometries. It is interesting to note that the proof of this classification
result does not involve the unitary invariants of the $n$-tuples obtained
earlier.
While the possible unitary invariants can be explicitly calculated,
at least for $n=2$ and $n=3$, the corresponding isometries are difficult
to identify.

\section{Model Multi-Isometries}

We begin with a few general remarks about sets of commuting isometries.
Consider a Hilbert space $\mathfrak{H}$, and a commutative semigroup
$\mathfrak{S}$ of isometric operators on $\mathfrak{H}$; i.e.,
$VW=WV\in\mathfrak{S}$
whenever $V,W\in\mathfrak{S}$. As in the Wold-von Neumann decomposition,
it was noted in \cite{key-6} that the subspace \[
\mathfrak{H}_{u}=\bigcap_{V\in\mathfrak{S}}V\mathfrak{H}\]
 is a reducing subspace for all the isometries in $\mathfrak{S}$ on which the
restrictions
 to it are unitary.
The semigroup $\mathfrak{S}$ is said to be \emph{completely nonunitary}
(or cnu) if $\mathfrak{H}_{u}=\{0\}$. Observe that individual elements
of $\mathfrak{S}$ might have unitary parts even if $\mathfrak{S}$
is cnu.

If the semigroup $\mathfrak{S}$ is generated by a set $\mathfrak{G}$ of
commuting
isometries, we will also call $\mathfrak{H}_{u}$ the unitary part
of $\mathfrak{G}$, and we will say that $\mathfrak{G}$ is cnu if
$\mathfrak{S}$ is cnu. If $\mathfrak{G}=\{ V_{1},V_{2},\dots,V_{n}\}$
is a finite set, observe that \[
\mathfrak{H}_{u}=\bigcap_{k=1}^{\infty}\bigcap_{i_{1},i_{2},\dots,i_{k}=1}^{n}V_
{i_{1}}V_{i_{2}}\cdots
V_{i_{k}}\mathfrak{H}=\bigcap_{k=1}^{\infty}V^{k}\mathfrak{H},\]
 where $V=V_{1}V_{2}\cdots V_{n}$; indeed, this is seen from the
inclusion \[
V_{i_{1}}V_{i_{2}}\cdots V_{i_{k}}\mathfrak{H}\subset V^{k}\mathfrak{H}\]
 Thus, as remarked earlier in \cite{key-1}, $\{ V_{1},V_{2},\dots,V_{n}\}$
is cnu if and only if $V$ is a unilateral shift (of arbitrary multiplicity).

For easier reference, an ordered $n$-tuple $(V_{1},V_{2},\dots,V_{n})$
of commuting isometries will be called simply an $n$-isometry. In
this paper we will focus mostly on completely nonunitary $n$-isometries. However, when the $n$-isometry $\{V_1,V_2,\ldots, V_n\}$ is not completely non-unitary, then the system of restrictions $V_i|_{{\cl H}_u}$ of $V_i$, $i=1,2,\ldots, n$,  will be called the unitary part of the $n$-isometry, while the system of restrictions to ${\cl H}\ominus {\cl H}_u$ will be called the cnu part.

We recall now that every shift can be realized conveniently as
an operator on a Hardy space. Thus, given a Hilbert space $\mathfrak{E}$,
consider the Hardy space $H^{2}(\mathfrak{E})$ of
Taylor series with square summable coefficients in $\mathfrak{E}$. The operator
$S_{\mathfrak{E}}$
defined by
\[
(S_{\mathfrak{E}}f)(z)=zf(z),\quad f\in H^{2}(\mathfrak{E}),\quad z\in {\bb D} = \{\lambda\in {\bb C}\colon \ |\lambda|<1\},
\]
 is a unilateral shift of multiplicity equal to the dimension of $\mathfrak{E}$.
It will be convenient to identify $\mathfrak{E}$ with the collection
of constant functions in $H^{2}(\mathfrak{E})$.

In view of the above considerations, it is of interest to study
isometries $V$ for which $S_{\mathfrak{E}}=VW$, where $W$ is some
other isometry commuting with $V$. We will call such an isometry
$V$ an isometric divisor of $S_{\mathfrak{E}}$. The following description
of divisors was first proved in \cite{key-1}. Our proof is somewhat
simpler.

For $P$ a projection, we set $P^\bot = I-P$.

\begin{thm}
The following are equivalent.

\emph{(1)} $V$ is an isometric divisor of $S_{\mathfrak{E}}$.

\emph{(2)} There exist a unitary operator $U$ on $\mathfrak{E}$
and an orthogonal projection $P$ on $\mathfrak{E}$ such that
$(Vf)(z)=U(zP+P^\bot)f(z)$
for all $f\in H^{2}(\mathfrak{E})$.
\end{thm}
\begin{proof}
It is trivial to verify that an operator $V$, as given in (2), is a divisor
of $S_{\mathfrak{E}}$ and, in fact, $VW=WV=S_{\mathfrak{E}}$, with
$(Wf)(z)=(P+zP^\bot)U^{*}f(z)$. Conversely, assume $W$ is an isometry
and $VW=WV=S_{\mathfrak{E}}$. Since $V$ and $W$ commute with $S_{\mathfrak{E}}$,
there exist inner $\mathfrak{L}(\mathfrak{E})$-valued functions
$\Theta,\Omega$ such that
$\Theta(z)\Omega(z)=\Omega(z)\Theta(z)=zI_{\mathfrak{E}}$,
$(Vf)(z)=\Theta(z)f(z)$, and $(Wf)(z)=\Omega(z)f(z)$ for $f\in
H^{2}(\mathfrak{E})$
and $|z|<1$. Now, the range of $W$ contains the range of $S_{\mathfrak{E}}$,
so  there exists a projection $P$ on $\mathfrak{E}$ for which
\[
WH^{2}(\mathfrak{E})=S_{\mathfrak{E}}H^{2}(\mathfrak{E})\oplus
P\mathfrak{E}=\Omega_{1}H^{2}(\mathfrak{E}),\]
 where $\Omega_{1}(z)=zP^\bot+P$. The Beurling-Lax-Halmos theorem
provides a unitary operator $U$ on $\mathfrak{E}$ such that
$\Omega_{1}(z)=\Omega(z)U$
for $|z|<1$. It is easy now to conclude for almost every $z$ on
the unit circle that \[
\Theta(z)=z\Omega(z)^{*}=zU\Omega_{1}(z)^{*}=U(P^\bot+zP),\]
 so that $V$ satisfies (2).
\end{proof}
We will denote by $V_{U,P}$ the divisor of $S_{\mathfrak{E}}$ described
in condition (2) of the above result.

\begin{lem}
Consider unitary operators $U,U_{1},U_{2}$, and orthogonal projections
$P,P_{1},P_{2}$ on the Hilbert space $\mathfrak{E}$. The following
are equivalent.

\emph{(1)} $V_{U,P}=V_{U_{1},P_{1}}V_{U_{2},P_{2}}$;

$(2)$ $U=U_{1}U_{2}$, and $P=P_{2}+U_{2}^{*}P_{1}U_{2}$.
\end{lem}
\begin{proof}
We only argue  the nontrivial implication (1)$\Rightarrow$(2).
Identification of coefficients yields the equations $U_{1}P_{1}U_{2}P_{2}=0$,
$U_{1}[P_{1}U_{2}P^\bot_{2}+P^\bot_{1}U_{2}P_{2}]=UP$, and
$U_{1}P^\bot_{1}U_{2}P^\bot_{2}=UP^\bot$.
The first two relations yield $U_{1}P_{1}U_{2}+U_{1}U_{2}P_{2}=UP$,
while the first and third yield
$U_{1}P^\bot_{1}U_{2}-V_{1}V_{2}P_{2}=VP^\bot$, where $V=V_{U,P}$.
Adding these last two equalities we obtain $U=U_{1}U_{2}$, and therefore
$U_{1}U_{2}P=VP=U_{1}P_{1}U_{2}+U_{1}U_{2}P_{2}$. This gives immediately
the second equality in (2).
\end{proof}
It should be noted that the relation $P=P_{2}+U_{2}^{*}P_{1}U_{2}$
contains implicitly the fact that the two projections on the right
hand side are orthogonal or, equivalently, $P_{1}U_{2}P_{2}=0$. Indeed,
the sum $Q=Q_{1}+Q_{2}$ of two orthogonal projections $Q_{1},Q_{2}$
is a contraction if and only if the ranges of $Q_{1}$ and $Q_{2}$
are orthogonal, and in this latter case $Q$ is the orthogonal projection
onto the span of the two ranges.

We can now give a model for arbitrary cnu $n$-isometries. Consider
a Hilbert space $\mathfrak{E}$, unitary operators $U_{1},U_{2},\dots,U_{n}$
on $\mathfrak{E}$, and orthogonal projections $P_{1},P_{2},\dots,P_{n}$
on $\mathfrak{E}$. The $n$-tuple
$(V_{U_{1},P_{1}},V_{U_{2},P_{2}},\dots,V_{U_{n},P_{n}})$
is called a \emph{model $n$-isometry} if the following conditions
are satisfied:

\begin{itemize}
\item[(a)] $U_{i}U_{j}=U_{j}U_{i}$ for $i,j=1,2,\dots,n$;
\item[(b)] $U_{1}U_{2}\cdots U_{n}=I_{\mathfrak{E}}$;
\item[(c)] $P_{j}+U_{j}^{*}P_{i}U_{j}=P_{i}+U_{i}^{*}P_{j}U_{i}\le I_{\mathfrak{E}}$
for $i\ne j$; and
\item[(d)]
$P_{1}+U_{1}^{*}P_{2}U_{1}+U_{1}^{*}U_{2}^{*}P_{3}U_{2}U_{1}+\cdots+U_{1}^{*}U_{
2}^{*}\cdots U_{n-1}^{*}P_{n}U_{n-1}\cdots U_{2}U_{1}=I_{\mathfrak{E}}$.
\end{itemize}
Observe again that the projections appearing in (d) must be pairwise
orthogonal. It follows from an inductive application of the preceding
lemma that a model $n$-isometry is indeed an $n$-isometry, and \[
V_{U_{1},P_{1}}V_{U_{2},P_{2}}\cdots V_{U_{n},P_{n}}=S_{\mathfrak{E}}.\]
 When $n=2$ we have $U_{2}=U_{1}^{*}$ and $P_{2}=I-U_{1}P_{1}U_{1}^{*}$,
so a model two-isometry is determined by $U_{1}$ and $P_{1}$. More
generally, a model $n$-isometry is determined by $U_{j},P_{j}$ for
$1\le j\le n-1$. The relevant conditions on these operators are as follows.

\begin{thm}
Let $\mathfrak{E}$ be a Hilbert space,  $(U_{j})_{j=1}^{n-1}$
 unitary operators on $\mathfrak{E}$, and $(P_{j})_{j=1}^{n-1}$
orthogonal projections on $\mathfrak{E}$. Assume that the following
conditions are satisfied:

\emph{(1)} $U_{i}U_{j}=U_{j}U_{i}$ for $i,j=1,2,\dots,n-1$;

\emph{(2)} $P_{j}+U_{j}^{*}P_{i}U_{j}=P_{i}+U_{i}^{*}P_{j}U_{i}\le
I_{\mathfrak{E}}$
for $i\ne j$, $1\le i,j\le n-1$; and

\emph{(3)} $P_{1}+U_{1}^{*}P_{2}U_{1}+\cdots+U_{1}^{*}U_{2}^{*}\cdots
U_{n-2}^{*}P_{n-1}U_{n-2}\cdots U_{2}U_{1}\le I_{\mathfrak{E}}$.

Then there exists a unique unitary operator $U_{n}$ on $\mathfrak{E}$
and a unique projection $P_{n}$ on $\mathfrak{E}$, such that
$(V_{U_{j},P_{j}})_{j=1}^{n}$
is a model $n$-isometry.
\end{thm}
\begin{proof}
Observe that the projections in condition (3) must be pairwise orthogonal.
An inductive application of the above lemma shows that the product
\[
V_{U_{1},P_{1}}V_{U_{2},P_{2}}\cdots V_{U_{n-1},P_{n-1}}\]
 is equal to $V_{U,P}$, where $U=U_{1}U_{2}\cdots U_{n-1}$, and
$P$ is the sum in the left hand side of condition (3). We can then
choose $U_{n}=U^{*}$ and $P_{n}=I-UPU^{*}$ to obtain a model $n$-isometry, and
this is clearly
the only possible choice.
\end{proof}

The conditions in this proposition are vacuously satisfied if $n=2$,
but are fairly stringent for larger $n$. We illustrate this in case
$n=3$, when the only conditions on $(U_{1},U_{2},P_{1},P_{2})$ are
that $U_{1}U_{2}=U_{2}U_{1}$ and\[
P_{2}+U_{2}^{*}P_{1}U_{2}=P_{1}+U_{1}^{*}P_{2}U_{1}\le I_{\mathfrak{E}}.\]
Observe that the data $(U_{1},U_{2},P_{1},P_{2})$ can be recovered
from $(U_{1},U_{2},P_{1},Q_{1})$, where $Q_{1}=U_{1}^{*}P_{2}U_{1}$,
and the required conditions can be written more easily in terms of the
mutually orthogonal projections $P_{1}$ and $Q_{1}$. Moreover, as will be seen below, the range of $Q_1$ is an invariant subspace for a suitably defined compression of $V_{U_1,P_1}$.

\begin{prop}
Consider unitary operators $U_{1},U_{2}$, and orthogonal projections
$P_{1},Q_{1}$ on a Hilbert space $\mathfrak{E}$ such that
$U_{1}U_{2}=U_{2}U_{1}$
and $P_{1}+Q_{1}\le I_{\mathfrak{E}}$. Let us also set \[
P_{2}=U_{1}Q_{1}U_{1}^{*},Q_{2}=U_{2}^{*}P_{1}U_{2},U=U_{1}U_{2},\:\text{and}\:
\, T_{1}=P^\bot_{1}U_{1}|P^\bot_{1}\mathfrak{E}.\]

\begin{enumerate}
\item We have $P_{2}\le P_{1}+Q_{1}$ if and only if $Q_{1}\mathfrak{E}$
is an invariant subspace for $T_{1}$.
\item We have $Q_{2}\le P_{1}+Q_{1}$ if and only if $Q_{1}\mathfrak{E}$
contains $P^\bot_{1}U_{2}^{*}P_{1}\mathfrak{E}$.
\item We have $P_{2}Q_{2}=0$ if and only if $Q_{1}\mathfrak{E}$ is contained
in $U^{*}P^\bot_{1}\mathfrak{E}$.
\end{enumerate}
In summary, the inequalities\[
P_{2}+U_{2}^{*}P_{1}U_{2}\le P_{1}+U_{1}^{*}P_{2}U_{1}\le I_{\mathfrak{E}}\]
hold if and only if $Q_{1}\mathfrak{E}$ is an invariant subspace
for $T_{1}$ such that\[
P^\bot_{1}U_{2}^{*}P_{1}\mathfrak{E}\subset Q_{1}\mathfrak{E}\subset
U^{*}P^\bot_{1}\mathfrak{E}.\]
In particular, if these inclusions hold, we have \[
P_{1}U[P^\bot_{1}U_{1}]^{n}P^\bot_{1}U_{2}^{*}P_{1}=0\:\text{for}\: n\ge0.\]

\end{prop}
\begin{proof}
Note that $P_{2}\le P_{1}+Q_{1}$ is equivalent to the inclusion
$P^\bot_{1}P_{2}\mathfrak{E}\subset Q_{1}\mathfrak{E}.$
Now (1) follows immediately since
$P^\bot_{1}P_{2}\mathfrak{E}=P^\bot_{1}U_{1}Q_{1}\mathfrak{E}$.
Similarly, $Q_{2}\le P_{1}+Q_{1}$ is equivalent to
$P^\bot_{1}Q_{2}\mathfrak{E}\subset Q_{1}\mathfrak{E}$,
which renders (2) obvious. Finally, $P_{2}Q_{2}=0$ if and only if
$P_{2}\mathfrak{E}\subset Q^\bot_{2}\mathfrak{E}$ or, equivalently,\[
Q_{1}\mathfrak{E}=U_{1}^{*}P_{2}\mathfrak{E}\subset
U_{1}^{*}Q^\bot_{2}\mathfrak{E}.\]
The inclusion in (3) follows because\[
U_{1}^{*}Q^\bot_{2}\mathfrak{E}=U_{1}^{*}Q^\bot_{2}U_{1}\mathfrak{E}=(I-U_{1
}^{*}Q_{2}U_{1})\mathfrak{E}=(I-U^{*}P_{1}U)\mathfrak{E}=U^{*}P^\bot_{1}
\mathfrak{E}.\]
For the last assertion in the statement, notice first that the condition
$P_{2}+Q_{2}\le I_{\mathfrak{E}}$ implies $P_{2}Q_{2}=0.$ Therefore,
this statement simply says that the invariant subspace for $T_{1}$,
generated by $P^\bot_{1}U_{2}^{*}P_{1}\mathfrak{E}$, is contained
in $U^{*}P^\bot_{1}\mathfrak{E}$.
\end{proof}

The preceding result shows the importance of the operator $T_1$ which we call the ``povital operator'' and consider further in the following section. Moreover, the result suggests introducing the spaces\[
\mathfrak{Q}_{1,\text{min}}=\bigvee_{k=0}^{\infty}T_{1}^{k}P^\bot_{1}U_{2}^{*}
P_{1}\mathfrak{E},\]
and \[
\mathfrak{Q}_{1,\text{max}}=P^\bot_{1}\mathfrak{E}\ominus\bigvee_{k=0}^{\infty
}T_{1}^{*k}P^\bot_{1}U_{1}^{*}U_{2}^{*}P_{1}\mathfrak{E}.\]
In other words, $\mathfrak{Q}_{1,\text{min}}$ is the smallest invariant
subspace for $T_{1}$ containing $P^\bot_{1}U_{2}^{*}P_{2}\mathfrak{E}$,
while $\mathfrak{Q}_{1,\text{max}}$ is the largest invariant subspace
for $T_{1}$ contained in $P^\bot_{1}\mathfrak{E}\ominus P^\bot_1U^*P_1\mathfrak{E}$.

\begin{thm}
Assume that $\mathfrak{E}$ is a Hilbert space, $U_{1},U_{2}$ are
commuting unitary operators on $\mathfrak{E}$ and $P_{1}$ is an
orthogonal projection on $\mathfrak{E}$. There exists a projection
$P_{2}$ on $\mathfrak{E}$ such that\[
P_{2}+U_{2}^{*}P_{1}U_{2}\le P_{1}+U_{1}^{*}P_{2}U_{1}\le I_{\mathfrak{E}}\]
if and only if
\emph{$\mathfrak{Q}_{1,\text{min}}\subset\mathfrak{Q}_{1,\text{max}}$.}
When this condition is satisfied, the general form of such projections
$P_{2}$ is $P_{2}=U_{1}Q_{1}U_{1}^{*}$, where $Q_{1}$ is the orthogonal
projection onto an invariant subspace of $T_{1}$ satisfying
\emph{
\[
\mathfrak{Q}_{1,\text{min}}\subset
Q_{1}\mathfrak{E}\subset\mathfrak{Q}_{1,\text{max}}.\]
}
\end{thm}
\begin{proof}
This is just a summary of the preceding discussion.
\end{proof}
When $\mathfrak{E}$ is finite dimensional, the statement can be made
more precise.

\begin{cor}
Under the hypotheses of the preceding theorem, assume that $\mathfrak{E}$
is finite dimensional and $\mathfrak{Q_{1,\text{{\rm min}
}}\subset\mathfrak{Q_{1,\text{{\rm max}}}}}$.
\begin{enumerate}
\item If $Q_{1}$ is an arbitrary projection onto an invariant subspace
of $T_{1}$, and
\emph{
\[
\mathfrak{Q}_{1,\text{min}}\subset
Q_{1}\mathfrak{E}\subset\mathfrak{Q}_{1,\text{max}},\]
}
then
\[
P_{2}+U_{2}^{*}P_{1}U_{2}=P_{1}+U_{1}^{*}P_{2}U_{1}\le I_{\mathfrak{E}},\]
where $P_{2}=U_{1}Q_{1}U_{1}^{*}.$ Thus, there exists a 3-isometry
of the form \[
(V_{U_{1},P_{1}},V_{U_{2},P_{2}},V_{U_{3},P_{3}})\]
such that $V_{U_{1},P_{1}}V_{U_{2},P_{2}}V_{U_{3},P_{3}}=S_{\mathfrak{E}}.$
\item The space \emph{$\mathfrak{Q}_{1,\text{min}}$} is $\{0\}$ if and
only if among the 3-isometries in \emph{(1)} there is one where
$V_{U_{2},P_{2}}$
is unitary.
\item The space \emph{$\mathfrak{Q}_{1,\text{max}}$} equals
$P^\bot_{1}\mathfrak{E}$
if and only if among the 3-isometries in \emph{(1)} there is one such
that $V_{U_{3},P_{3}}$ is unitary.
\end{enumerate}
\end{cor}
\begin{proof}
The last part of the first statement can be deduced from the fact that the
projections
$P_{2}+U_{2}^{*}P_{1}U_{2}$ and $P_{1}+U_{1}^{*}P_{2}U_{1}$ have
the same rank, equal to the sum of the ranks of $P_{1}$ and $P_{2}$.
The second statement follows from the fact that $Q_{1}$ can be chosen
to be zero if and only if $\mathfrak{Q}_{1,\text{min}}=\{0\}$. Likewise,
the third statement follows from the fact that $Q_{1}$ can be chosen
to be $P^\bot_{1}$ if and only if
$\mathfrak{Q}_{1,\text{max}}=P^\bot_{1}\mathfrak{E}$.
\end{proof}
One could ask whether part (1) of the preceding corollary is true
when the dimension of $\mathfrak{E}$ is infinite. Unfortunately,
the answer is negative, as shown by the following example.

\begin{example}
Denote by $\mathfrak{E}=L^{2}$ the usual space of square integrable
functions on the unit circle (relative to normalized arclength measure),
and define unitary operators $U_{1},U_{2}$ on $\mathfrak{E}$ by
setting\[
(U_{1}f)(z)=\varphi(z)f(z),(U_{2}f)(z)=zf(z),\quad
f\in\mathfrak{E},z\in\partial\mathbb{D},\]
where the function $\varphi$ is defined to be equal to 1 on the upper
half-circle, and $-1$ on the lower half-circle. Consider also the
orthogonal projection $P_{1}$ on $\mathfrak{E}$ such that
$P_{1}\mathfrak{E}=(H^{2})^{\perp}.$
We claim that with these choices we have \[
\mathfrak{Q}_{1,\text{min}}=\mathfrak{Q}_{1,\text{max}}=\{0\},\]
and for the (unique) choice $Q_{1}=0$ we have\[
P_{2}+U_{2}^{*}P_{1}U_{2}\lvertneqq P_{1}+U_{1}^{*}P_{2}U_{2}.\]

Indeed, $U_{2}$ leaves $H^{2}\subset\mathfrak{E}$ invariant. Therefore, $(P^\bot_{1})U_{2}^{*}P_{1}=0$,
which in turn implies that $\mathfrak{Q}_{1,\text{min}}=\{0\}$. On
the other hand,
\[
(P^\bot_{1})\mathfrak{E}\ominus\mathfrak{Q}_{1,\text{max}}=\bigvee_{n=0}^{\infty}
T_{1}^{*n} (P^\bot_{1})U_{1}^{*}U_{2}^{*}P_{1}\mathfrak{E}\]
 is the smallest invariant subspace for the operator
$T_{1}^{*}=T_{1}=T_{\varphi}$
generated by the space
$(P^\bot_{1})U_{1}^{*}U_{2}^{*}(H_{2})^{\perp}=P_{H^{2}}[\overline{\varphi
z}(H^{2})^{\perp}].$
This space is actually dense in $H^{2}$ already. Indeed, consider
a function $u\in H^{2}\ominus P_{H^{2}}[\overline{\varphi
z}(H^{2})^{\perp}]=H^{2}\cap(\varphi zH^{2}).$
There must exist $v\in H^{2}$ such that $u=\varphi zv$. The F. and
M. Riesz theorem implies (by looking at the upper half-circle) that
$u=zv$ and (looking at the lower half-circle) $u=-zy$, and therefore
$u=0$. It is now easy to see that\[
P_{2}+U_{2}^{*}P_{1}U_{2}=U_{2}^{*}P_{1}U_{2}=P_{\overline{z}(H^{2})^{\perp}}\ne
P_{1}+U_{1}^{*}P_{2}U_{1}=P_{1}=P_{(H^{2})^{\perp}}.\]

\end{example}
The discussion above shows that constructing cnu 3-isometries can
be a delicate task. If
$\mathfrak{Q}_{1,\text{max}}\supset\mathfrak{Q}_{1,\text{min}}$
and the compression of the pivotal operator $T_{1}$ in Theorem 2.5 to the
space $\mathfrak{Q}_{1,\text{max}}\ominus\mathfrak{Q}_{1,\text{min}}$
is transitive (e.g., when this space has dimension zero or one), then
the only choices for $Q_{1}$ are the projections onto
$\mathfrak{Q}_{1,\text{max}},\mathfrak{Q}_{1,\text{min}}$,
and even these may fail to produce 3-isometries, as seen in the preceding
example. It is, however, possible to formulate equivalent conditions
for the existence of $P_{2}$, given $U_{1},U_{2}$ and $P_{1}$.
To find these conditions assume that, given this data, the inclusion
$\mathfrak{Q_{1,\text{{max}}}\supset\mathfrak{Q_{1,\text{{min}}}}}$ is satisfied.
According to Theorem 2.5, the inequalities
\[
P_{2}+U_{2}^{*}P_{1}U_{2}\le P_{1}+U_{1}^{*}P_{2}U_{1}\le I_{\mathfrak{E}}\]
are satisfied when $P_{2}=U_{1}Q_{1}U_{1}^{*}$, and $Q_{1}$ is the
orthogonal projection on either of the spaces
$\mathfrak{Q}_{1,\text{max}},\mathfrak{Q}_{1,\text{min}}$.
This allows us to define an isometric operator\[
W:P_{1}+\mathfrak{Q}_{1,\text{max}}\to P_{1}+\mathfrak{Q}_{1,\text{max}}\]
by setting \[
W(x_{1}+x_{2})=U_{2}^{*}x_{1}+U_{1}x_{2},\quad x_{1}\in
P_{1}\mathfrak{{E}},x_{2}\in\mathfrak{Q_{1,\text{max}}}.\]

\begin{thm}
Assume that $\mathfrak{E}$ is a Hilbert space, $U_{1},U_{2}$ are commuting
unitary operators on $\mathfrak{E}$, $P_{1}$ is an orthogonal projection,
and $\mathfrak{Q}_{1,\text{{\rm min} }}\subset\mathfrak{Q}_{1,\text{{\rm max}}}$.
Define the isometry $W$ as above and consider the operator $T_{1}$
used in Proposition 2.4.
\begin{enumerate}
\item A subspace $\mathfrak{Q_{1}}$ such that $\mathfrak{Q}_{1,\text{{\rm min}
}}\subset\mathfrak{Q}_{1}\subset\mathfrak{Q}_{1,\text{{\rm max}}}$
is invariant for $T_{1}$ if and only if $P_{1}\mathfrak{E}+\mathfrak{Q}_{1}$
is invariant for $W$.
\item Let $\mathfrak{Q}_{1}$ be an invariant subspace of $T_{1}$ such that
$\mathfrak{Q}_{1,\text{{\rm min} }}\subset\mathfrak{Q}_{1}\subset\mathfrak{Q}_{1,\text{{\rm
max}}}$,
denote by $Q_{1}$ the orthogonal projection onto $\mathfrak{Q}_{1}$,
and set $P_{2}=U_{1}Q_{1}U_{1}^{*}$. We have\[
P_{2}+U_{2}^{*}P_{1}U_{2}=P_{1}+U_{1}^{*}P_{2}U_{1}\le I_{\mathfrak{E}}\]
if and only if $W|(P_{1}\mathfrak{E}+\mathfrak{Q}_{1})$ is a unitary operator.
\item There exists an orthogonal projection $P_{2}$ on $\mathfrak{E}$ such
that\[
P_{2}+U_{2}^{*}P_{1}U_{2}=P_{1}+U_{1}^{*}P_{2}U_{1}\le I_{\mathfrak{E}}\]
 if and only if $P_{1}\mathfrak{E}+\mathfrak{Q}_{1,\text{{\rm min}}}$ is
contained in the  unitary part of $W$ in the von Neumann---Wold decomposition.
The collection of such projections $P_{2}$ is a complete lattice.
\end{enumerate}
\end{thm}
\begin{proof}
Assume first that $\mathfrak{Q}_{1}$ is invariant for $T_{1}$. As noted
above, $W$ leaves $P_{1}\mathfrak{E}+\mathfrak{Q}_{1,\text{min}}$ invariant.
Therefore, for $p_{1}\in P_{1}\mathfrak{E}$ and $q_{1}\in\mathfrak{Q}_{1}$
we have \[
W(p_{1}+q_{1})=Wp_{1}+P_{1}U_{1}q_{1}+T_{1}q_{1}\in
P_{1}\mathfrak{E}+\mathfrak{Q}_{1,\text{min}}+\mathfrak{Q}_{1}\subset
P_{1}\mathfrak{E}+\mathfrak{Q}_{1}.\]
Conversely, if $W$ leaves $P_{1}\mathfrak{E}+\mathfrak{Q}_{1}$ invariant,
the above formula can be rewritten as\[
T_{1}q_{1}=W(p_{1}+q_{1})-Wp_{1}\in P_{1}\mathfrak{E}+\mathfrak{Q}_{1},\]
and this clearly implies that $T_{1}$ leaves $\mathfrak{Q}_{1}$ invariant.
This proves (1). To verify (2), one only needs to observe that the
range of $P_{2}+U_{2}^{*}P_{1}U_{2}$ is precisely equal to
$W(P_{1}\mathfrak{E}+\mathfrak{Q}_{1})$.
Part (3) follows immediately from (1) and (2).
\end{proof}
The preceding result also clarifies Example 2.7, for which the isometry
$W$ is pure.

The following result appears in a slightly different form as
Theorem 3.2 in \cite{key-1}.

\begin{thm}
${}$
\begin{enumerate}
\item Any cnu $n$-isometry is unitarily equivalent to a model $n$-isometry.
\item Consider two model $n$-isometries $(V_{U_{j},P_{j}})_{j=1}^{n}$
and $(V_{U'_{j},P'_{j}})_{j=1}^{n}$, where $U_{j},P_{j}$ act on
$\mathfrak{E}$ and $U'_{j},P'_{j}$ act on $\mathfrak{E}'$. These
model $n$-isometries are unitarily equivalent if and only if there
exists a unitary operator $W:\mathfrak{E}\to\mathfrak{E}'$ satisfying
$WU_{j}=U'_{j}W$ and $WP_{j}=P'_{j}W$ for all $j$.
\end{enumerate}
\end{thm}
\begin{proof}
Let $(V_{1},V_{2},\dots,V_{n})$ be a cnu $n$-isometry. Up to unitary
equivalence, we may assume that $V_{1}V_{2}\cdots V_{n}=S_{\mathfrak{E}}$
for some Hilbert space $\mathfrak{E}$. It follows then that each
$V_{j}$ is of the form $V_{U_{j},P_{j}}$, and properties (a)---(d)
follow from the Lemma 2.2. The second part of the statement follows from
the fact that any unitary operator $X:H^{2}(\mathfrak{E})\to
H^{2}(\mathfrak{E}')$
satisfying $XS_{\mathfrak{E}}=S_{\mathfrak{E}'}X$ must be the multiplication
operator defined by some (constant) unitary $W:\mathfrak{E}\to\mathfrak{E}'$.
\end{proof}
It is natural to characterize various properties of $n$-isometries
in terms of the corresponding models. This is thoroughly pursued in
\cite{key-8}. Here we only note the following result. We recall that
operators $T$ and $S$ are said to doubly commute if $TS=ST$ and
$T^{*}S=ST^{*}$.

\begin{prop}
A model two-isometry $(V_{U_{1},P_{1}},V_{U_{2},P_{2}})$ consists
of doubly commuting operators if and only if
$P_{1}U_{1}^{*}P_{1}=U_{1}^{*}P_{1}$,
i.e., $U_{1}^{*}$ leaves the range of $P_{1}$ invariant. Furthermore, these conditions are satisfied if and only the pivotal operator $T_1$ is an isometry.
\end{prop}
\begin{proof}
Fix $f\in H^{2}(\mathfrak{E})$ and $|z|<1$. A calculation shows
that \[
(V_{U_{2},P_{2}}^{*}f)(z)=U_{1}P^\bot_{1}\frac{f(z)-f(0)}{z}+U_{1}P_{1}f(z),\]
 and further computation yields \[
[(V_{U_{2},P_{2}}^{*}V_{U_{1},P_{1}}-V_{U_{1},P_{1}}V_{U_{2},P_{2}}^{*})f](z)=U_
{1}P_{1}U_{1}P^\bot_{1}f(0).\]
 Thus double commutation is equivalent to $P_{1}U_{1}P^\bot_{1}=0$.
\end{proof}
A related result is observed in \cite{key-1}, namely that the commutator
$V_{1}^{*}V_{2}-V_{2}V_{1}^{*}$ is compact if and only if
$P_{1}U_{1}P^\bot_{1}$
is compact.

The preceding result extends in the obvious way to arbitrary model
$n$-isometries. The condition for double commutativity is simply
that, for each $j$, $U_{j}^{*}$ leaves the range of $P_{j}$ invariant.

\section{The Invariants of Bi-Isometries}

As noted above, a complete set of unitary invariants of cnu bi-isometries
is provided by triples $(\mathfrak{E},U,P)$, where $\mathfrak{E}$
is a Hilbert space and $U$ and $P$ are operators on $\mathfrak{E}$
with $U$ unitary and $P$ an orthogonal projection. The model operators are $V_{U_1,P_1}$ and $V_{U_2,P_2}$, where $U_1 = U$, $P_1=P$, $U_2=U^*$, and $P_1 = UP^\bot U^*$. For easier reference,
we will call such a triple a model triple. Two model triples
$(\mathfrak{E},U,P)$
and $(\mathfrak{E}_{1},U_{1},P_{1})$ determine unitarily equivalent
bi-isometries if and only if there exists a unitary operator
$A:\mathfrak{E}\to\mathfrak{E}_{1}$
satisfying $AU=U_{1}A$ and $AP=P_{1}A$.

Fix a model triple $(\mathfrak{E},U,P)$, and introduce auxiliary
spaces $\mathfrak{F}=P^\bot\mathfrak{E}$, $\mathfrak{D}=(\mathfrak{F}\vee
U\mathfrak{F})\ominus\mathfrak{F}$,
$\mathfrak{F}'=\mathfrak{E}\ominus(\mathfrak{F}\vee U\mathfrak{F})$,
$\mathfrak{D}_{*}=(\mathfrak{F}\vee U^{*}\mathfrak{F})\ominus\mathfrak{F}$,
and $\mathfrak{F}'_{*}=\mathfrak{E}\ominus(\mathfrak{F}\vee U^{*}\mathfrak{F})$.
Observe that $\mathfrak{F}\vee U\mathfrak{F}=\mathfrak{F}\oplus\mathfrak{D}$,
$\mathfrak{F}\vee U^{*}\mathfrak{F}=\mathfrak{F}\oplus\mathfrak{D}_{*}$,
$\mathfrak{D}\oplus\mathfrak{F}'=\mathfrak{D}_{*}\oplus\mathfrak{F}'_{*}$,
$U(\mathfrak{F}\vee U^{*}\mathfrak{F})=\mathfrak{F}\vee U\mathfrak{F}$,
and consequently $U\mathfrak{F}'_{*}=\mathfrak{F}$. Consider also
the contraction operator $T$ on $\mathfrak{F}$ defined by $T=P^\bot U|\mathfrak{F}$.
(Note that $T$ coincides with the pivotal operator $T_1$ in Proposition 2.4.) Our goal in this section is to understand better the relation of $T$ to the
corresponding bi-isometry.

It is well-known (cf. for example Theorem IV.3.1 in \cite{ff}) that
the unitary operator
$U|\mathfrak{F}\oplus\mathfrak{D}_{*}:\mathfrak{F}\oplus\mathfrak{D}_{*}\to
\mathfrak{F}\oplus\mathfrak{D}$
is essentially given by the Julia-Halmos matrix \[
\left[\begin{array}{cc}
T & D_{T^{*}}\\
D_{T} & -T^{*}\end{array}\right],\]
where $D_{T}=(I-T^{*}T)^{1/2}$ and $D_{T^{*}}=(I-TT^{*})^{1/2}$.
More precisely, there are unique unitary operators
$W:\mathfrak{D}_{T}=(D_{T}\mathfrak{E})^{-}\to\mathfrak{D}$
and $W_{*}:\mathfrak{D}_{T^{*}}\to\mathfrak{D}_{*}$ such that, for
$f\in\mathfrak{F}$ and $d_{*}\in\mathfrak{D}_{*}$ we have
\[
U(f\oplus d_{*})=[Tf+D_{T^*}W_{*}^{*}d_{*}]\oplus[WD_{T}f-WT^{*}W_{*}^{*}d_{*}].\]
 To simplify the notation, we  replace the original model triple by
the equivalent model triple
$(\mathfrak{F}\oplus\mathfrak{D}_{T}\oplus\mathfrak{F}',\Omega^{*}U\Omega,\Omega
^{*}P\Omega=P_{\mathfrak{F}\oplus\{0\}\oplus\{0\}})$,
where
$\Omega:\mathfrak{F}\oplus\mathfrak{D}_{T}\oplus\mathfrak{F}'\to\mathfrak{E}$
is given by $\Omega(f\oplus d\oplus f')=f+Wd+f'$. If we also consider
the unitary
$\Omega':\mathfrak{F}\oplus\mathfrak{D}_{T^{*}}\oplus\mathfrak{F}'\to\mathfrak{E
}$
given by $\Omega'(f\oplus d_{*}\oplus f')=f+W_*d_{*}+U^{*}f'$, the
new unitary $\Omega^{*}U\Omega$ can be factored as \[
\Omega^{*}U\Omega=(\Omega^{*}U\Omega')(\Omega^{\prime*}\Omega),\]
 and now the operator
$\Omega^{*}U\Omega^\prime:\mathfrak{F}\oplus\mathfrak{D}_{T^{*}}\oplus\mathfrak{F}'\to
\mathfrak{F}\oplus\mathfrak{D}_{T}\oplus\mathfrak{F}'$
is represented by the matrix
\[
\Omega^{*}U\Omega^\prime=\left[\begin{array}{ccc}
T & D_{T^{*}} & 0\\
D_{T} & -T^{*} & 0\\
0 & 0 & I_{\mathfrak{F}'}\end{array}\right],\]
 while
$\Omega^{\prime*}\Omega:\mathfrak{F}\oplus\mathfrak{D}_{T}\oplus\mathfrak{F}'\to
\mathfrak{F}\oplus\mathfrak{D}_{T^{*}}\oplus\mathfrak{F}'$
must have the form \[
\Omega^{\prime*}\Omega=\left[\begin{array}{cc}
I_{\mathfrak{F}} & 0\\
0 & Z\end{array}\right],\]
where
$Z:\mathfrak{D}_{T}\oplus\mathfrak{F}'\to\mathfrak{D}_{T^{*}}\oplus\mathfrak{F}'
$
is a unitary. We summarize this construction in the following result.

\begin{prop}
Consider Hilbert spaces $\mathfrak{F},\mathfrak{F}'$, a contraction
$T$ on $\mathfrak{F}$, and a unitary operator
$Z:\mathfrak{D}_{T}\oplus\mathfrak{F}'\to\mathfrak{D}_{T^{*}}\oplus\mathfrak{F}'
$.
Associated with this data is a model triple $(\mathfrak{E},U,P)$,
where $\mathfrak{E}=\mathfrak{F}\oplus\mathfrak{D}_{T}\oplus\mathfrak{F}'$,
$P^\bot\mathfrak{E}=\mathfrak{F}\oplus\{0\}\oplus\{0\}$, and $U=W_{1}W_{2}$,
with
$W_{1}:\mathfrak{F}\oplus\mathfrak{D}_{T^{*}}\oplus\mathfrak{F}'\to\mathfrak{F}
\oplus\mathfrak{D}_{T}\oplus\mathfrak{F}'$
and $W_{2}:\mathfrak{F}\oplus\mathfrak{D}_{T}\oplus\mathfrak{F}'\to\mathfrak{F}
\oplus\mathfrak{D}_{T^{*}}\oplus\mathfrak{F}'$
are given by the matrices \[
W_{1}=\left[\begin{array}{ccc}
T & D_{T^{*}} & 0\\
D_{T} & -T^{*} & 0\\
0 & 0 & I_{\mathfrak{F}'}\end{array}\right],\quad W_{2}=\left[\begin{array}{ccc}
I_{\mathfrak{F}} & 0&0\\
0\\
0&\multicolumn{2}{c}{\raisebox{1.5ex}[0cm][0cm]{\text{\LARGE $Z$}}} \end{array}\right].\]

\begin{enumerate}
\item Every model triple can be obtained, up to unitary equivalence, in
the manner described above.
\item The data $(\mathfrak{F},\mathfrak{F}',T,Z)$ and
$(\mathfrak{F}_{1},\mathfrak{F}_{1}',T_{1},Z_{1})$
determine equivalent model triples if and only if there exist unitary
operators $A:\mathfrak{F}\to\mathfrak{F}_{1}$ and
$A':\mathfrak{F}'\to\mathfrak{F}'_{1}$
satisfying $AT=T_{1}A$ and \[
((A\oplus A')|\mathfrak{D}_{T}\oplus\mathfrak{F}')Z=Z_{1}(A\oplus
A')|\mathfrak{D}_{T^{*}}\oplus\mathfrak{F}'\]
.
\end{enumerate}
\end{prop}
The uniqueness assertion in (2) follows from the fact that the construction
of $\mathfrak{F},\mathfrak{F}',T$, and $Z$ from $(\mathfrak{E},U,P)$
is invariant under unitary equivalence.

It is instructive to give the explicit form of the model operators $V_{U_1,P_1}$ and $V_{U_2,P_2}$ in terms of the operators $T$ and $Z$ in Proposition 3.1. So we will take
\[
\mathfrak{E} = \mathfrak{F} \oplus \mathfrak{D}_T \oplus \mathfrak{F}'
\]
and consequently also
\[
H^2(\mathfrak{E}) = H^2(\mathfrak{F}) \oplus H^2(\mathfrak{D}_T) \oplus H^2(\mathfrak{F}').
\]
The operator-valued functions appearing in the definition of $V_{U_1,P_1}$ and $V_{U_2,P_2}$ are
\[
U_1(zP_1 + P^\bot_1) = U(zP+P^\bot)
\]
and
\[
U_2(zP_2 + P^\bot_2) = (P+zP^\bot)U^*,
\]
respectively. In matrix form the first is
\begin{align*}
&\left[\begin{matrix}
T&D_{T^*}&0\\ D_T&-T^*&0\\ 0&0&I_{\mathfrak{F}'}\end{matrix}\right] \begin{bmatrix} I_{\mathfrak{F}}&0&0\\ 0\\ 0&\multicolumn{2}{c}{\raisebox{1.5ex}[0cm][0cm]{\text{\LARGE $Z$}}}\end{bmatrix} \begin{bmatrix} I_{\mathfrak{F}}&0&0\\ 0\\ 0&\multicolumn{2}{c}{ \raisebox{1.5ex}[0cm][0cm]{\text{\LARGE $zI_{\mathfrak{E}\ominus\mathfrak{F}}$}}}\end{bmatrix} =\\
&\quad = \begin{bmatrix} \mathfrak{T}&D_{\mathfrak{T}^*}&0\\ D_{\mathfrak{T}}&-T^*&0\\ 0&0&I_{\mathfrak{F}'}\end{bmatrix} \begin{bmatrix} I_{\mathfrak{F}}&0&0\\ 0\\ 0&\multicolumn{2}{c}{\raisebox{1.5ex}[0cm][0cm]{\text{\LARGE $zZ$}}}\end{bmatrix}
\end{align*}
and the second is
\[
\begin{bmatrix} zI_{\mathfrak{F}}&0&0\\ 0\\ 0&\multicolumn{2}{c}{\raisebox{1.5ex}[0cm][0cm]{\text{\LARGE $Z^*$}}}\end{bmatrix}
\begin{bmatrix} T^*&D_{\mathfrak{T}}&0\\ D_{\mathfrak{T}^*}&-T&0\\ 0&0&I_{\mathfrak{F}'}.
\end{bmatrix}.
\]
Consequently, we have
\[
\left\|V^*_2\begin{bmatrix} f\\ d\\ \varphi\end{bmatrix}\right\|^2 = \|g\|^2 + \left\|\begin{bmatrix} d\\ \varphi\end{bmatrix}\right\|^2,
\]
where $g\in H^2(\mathfrak{E})$ is given by
\[
g(z) = [f(z) - f(0)]/z\qquad (0\ne z\in {\bb D})
\]
and where
\[
f\in H^2(\mathfrak{E}), d\in H^2(\mathfrak{D}_T),\qquad \varphi\in H^2(\mathfrak{F}').
\]
Thus
\[
\text{ker } V^*_2 =\mathfrak{F},
\]
where $\mathfrak{F}$ \emph{is viewed as the subspace of $H^2(\mathfrak{E})$ formed by the constant functions with values in $\mathfrak{F}$}.  It readily follows that
\[
V^*_1|\text{ker } V^*_2 = T^*
\]
and hence
\[
T = P_{\text{ker } V^*_2} V_1|\text{ker } V^*_2.
\]
Thus we obtain the following result.

\begin{prop}
The family of all operators of the form $P_{\text{ker } V^*_2} V_1|\text{ker } V^*_2$, when $\{V_1,V_2\}$ runs over all c.n.u.\ bi-isometries, is (up to a unitary equivalence) the family of all contractions in Hilbert spaces.
\end{prop}

Continuing our study of the operators $V_1,V_2$, we next introduce the space $\mathfrak{F}_u$ of the unitary part of $T$ and notice that
\[
V_1\begin{bmatrix} f\\ 0\\ 0\end{bmatrix} = \begin{bmatrix}
Tf\\ 0\\ 0\end{bmatrix}\quad f\in H^2(\mathfrak{F}_u)
\]
and hence
\[
V_1(H^2(\mathfrak{F}_u) \oplus\{0\}\oplus \{0\}) = H^2(\mathfrak{F}) \oplus \{0\} \oplus \{0\}.
\]
Thus, $H^2(\mathfrak{F}_u)$ \emph{viewed as a subspace of} $H^2(\mathfrak{E})$, is the subspace $H^2(\mathfrak{E})^{(1)}_u$ of the unitary part of $V_1$ and
\[
(V_1f)(z) = Tf(z)\qquad z\in {\bb D}, f\in H^2(\mathfrak{F}_u).
\]
Therefore, $H^2(\mathfrak{F}_u)$ reduces $V_1$.

On the other hand we have
\[
V^n_2\begin{bmatrix}
f\\ 0\\ 0
     \end{bmatrix} = \begin{bmatrix}
z^nT^{*^n}f\\ 0\\ 0\end{bmatrix}\quad n=0,1,\ldots, f\in H^2(\mathfrak{F}_u).
\]
and consequently
\[
V_2H^2(\mathfrak{F}_u) \subset H^2(\mathfrak{F}_u).
\]
But
\[
\bigvee_{n\ge 0} V^n_2\mathfrak{F} = \mathfrak{F} \oplus V_2\mathfrak{F} \oplus V^2_2\mathfrak{F} \oplus\cdots
\]
is the space of the c.n.u.\ part of $V_2$ and so, since
\[
H^2(\mathfrak{F}_u) = \mathfrak{F}_u \oplus V_2\mathfrak{F}_u \oplus V^2_2 \mathfrak{F}_u \oplus\cdots,
\]
it is obvious that $H^2(\mathfrak{F}_u)$ reduces $V_2$. We have thus obtained the following.

\begin{lem}
Let $\mathfrak{F}_u$ be the space of the unitary part of $T$. Then $H^2(\mathfrak{F}_u)$ reduces both $V_1,V_2$ and is included in the space $H^2(\mathfrak{E})^{(1)}_u$ of the unitary part of $V_1$ and in that of the c.n.u.\ part of $V_2$, i.e.\ $H^2(\mathfrak{E})^{(2)}_{cnu}$.
\end{lem}

In our further investigation of the structure of the operators $V_1,V_2$, we can now restrict our attention to the restrictions of $V_1$ and $V_2$ to $H^2(\mathfrak{E})\ominus H^2(\mathfrak{F}_u)$. That means that, without loss of generality, we can assume that $\mathfrak{F}_u = \{0\}$ during this investigation. Then for $f\in \mathfrak{F}$ and $g\in H^2(\mathfrak{E})^{(1)}_u$ we have (with $\langle\cdot,\cdot\rangle$ denoting the scalar product in $H^2(\mathfrak{E})$)
\begin{align*}
\langle f,g\rangle &= \langle f,V^n_1V^{*^n}_1g\rangle = \langle V^{*^n}_1f, V^{*^n}_1g\rangle =\\
&= \langle T^{*^n}f, V^{*^n}_1g\rangle \to 0
\end{align*}
if $\|T^{*^n}f\|\to 0$ for $n\to\infty$.

We will consider now the case when the latter convergence holds for all $f\in \mathfrak{F}$, that is, the case when
\[
T\in C_{\cdot 0}.
\]
The above calculation shows  in this case that
\[
\mathfrak{F} \perp H^2(\mathfrak{E})^{(1)}_u.
\]
But
\[
H^2(\mathfrak{E})^{(2)}_{cnu} = \bigvee^\infty_{n=0} V^n_2\mathfrak{F} = \mathfrak{F} \oplus V_2 \mathfrak{F}\oplus V^2_2\mathfrak{F}\oplus\cdots
\]
and for $n\ge 1$
\[
\ker V^{*^n}_2 = \mathfrak{F} \oplus V_2\mathfrak{F} \oplus\cdots\oplus V^{n-1}_2\mathfrak{F}.
\]
Thus
\[
V^*_1\ker V^{*^n}_2 \subset \ker V^{*^n}_2
\]
and consequently
\[
V^*_1 H^2(\mathfrak{E})^{(2)}_{cnu} \subset H^2(\mathfrak{E})^{(2)}_{cnu}.
\]

At this stage in our study we need the following

\begin{lem}
Let $A$ on the Hilbert space $H^2(\mathfrak{G})$ be contractive commuting with the canonical shift $S_{\mathfrak{E}}$. If $A^{*^n}|\ker S^*_{\mathfrak{E}}\to 0$ strongly for $n\to \infty$, then $A^{*^n}\to 0$ strongly too.
\end{lem}

\begin{proof}
The operator $A$ is the multiplication operator on $H^2(\mathfrak{G})$ given by a bounded analytic operator-valued function
\[
A(z) = A_0+zA_1+\cdots,\quad \|A(z)\| \le 1 \quad (z\in {\bb D})
\]
on ${\bb D}$, where the $A_j$ are operators on $\mathfrak{G}$. The condition in the statement means that $A^{*^n}_0\to 0$ strongly. For $g(z) = g_0+zg_1 +\cdots z^Ng_N$, we have that $A^{*^n}g$ has the form
\[
(A^{*^n}g)(z) = g^{(n)}_0 + zg^{(n)}_1 +\cdots+ z^{N-1}g^{(n)}_{N-1} + z^NA^{*^n}_0g_N.
\]
Thus
\[
{\cl L} = \underset{n\to \infty}{\overline{\text{lim}}}\|A^{*^n}g\|^2 = \limsup\|A^{*^n}(g-S^N_{\mathfrak{E}}g_N)\|.
\]
Iterating this argument we obtain finally that
\[
{\cl L} = \lim\|A^{*^n}_0g_0\| = 0.\qquad \qed
\]
\renewcommand{\qed}{}\end{proof}

Returning to our study, we conclude that
\[
V^{*^n}_1|H^2(\mathfrak{E})^{(2)}_{cnu}\to 0 \text{ strongly}.
\]
The argument establishing the orthogonality $\mathfrak{F}\perp H^2(\mathfrak{E})^{(1)}_u$ now implies that
\[
H^2(\mathfrak{E})^{(2)}_{cnu} \perp H^2(\mathfrak{E})^{(1)}_u.
\]
Therefore
\[
H^2(\mathfrak{E})^{(2)}_u = H^2(\mathfrak{E}) \ominus H^2(\mathfrak{E})^{(2)}_u \supset H^2(\mathfrak{E})^{(1)}_u.
\]
But for any $h\in H^2(\mathfrak{E})^{(2)}_u$ and any $n=1,2,\ldots$ we have
\[
h = V^n_2h_n, \text{ where } h_n\in H^2(\mathfrak{E})^{(2)}_u
\]
so
\[
V^n_1h = (V_1V_2)^nh_n = z^nh_n \in z^nH^2(\mathfrak{E}).
\]
Therefore
\[
H^2(\mathfrak{E})^{(1)}_u = V^n_1H^2(\mathfrak{E})^{(1)}_u \subset z^nH^2(\mathfrak{E})
\]
for all $n=1,2,\ldots$~. This implies
\[
H^2(\mathfrak{E})^{(1)}_u = \{0\}.
\]

Returning to the case when $\mathfrak{F}_u$ may not be  $\{0\}$, we have thus established the following structure theorem.

\begin{thm}
Let $T = T_{cnu}\oplus T_u$ be the canonical decomposition of $T$ into its cnu part
 $T_{cnu}$ (on $\mathfrak{F}_{cnu}$) and unitary
 part $T_u$ (on $\mathfrak{F}_u$). In case
\[
T_{cnu}\in C_{\cdot 0},
\]
then $H^2(\mathfrak{F})(\subset H^2(\mathfrak{E}))$ is the space of
the unitary part of $V_1$.

This space is also reducing for $V_2$ and is contained  in the space
of the unitary part of $V_2$.
\end{thm}

Recalling that any bi-isometry $\{V_1,V_2\}$ is an orthogonal sum of its model\break $\{V_{U_1,P_1}, V_{U_1P_2}\}$ and a bi-isometry formed by unitary operators, we infer the following.

\begin{cor}
Let $\{V_1,V_2\}$ be any bi-isometry. If the c.n.u.\ part of
\[
V^*_1|\ker V^*_2
\]
is a $C_{10}$ contraction, then the Wold decomposition of $V_1$ reduces $V_2$ too.
\end{cor}

\begin{rem}
1)~~Let $\{V_1,V_2\}$  be a bi-isometry on $\mathfrak{H}$ such that the Wold decomposition of $V_1$ also reduces $V_2$. Let $\mathfrak{H} = \mathfrak{H}^{(1)}_{cnu} \oplus \mathfrak{H}^{(1)}_u$ be that decomposition. Thus
\[
\ker V^*_2 = \ker(V_2|\mathfrak{H}^{(1)}_{cnu})^* \oplus \ker(V_2|\mathfrak{H}^{(1)}_u)^*
\]
and (with a little abuse of notation)
\begin{align*}
T^* &= V^*_1|\ker V^*_2 = [(V^*_1|\mathfrak{H}^{(1)}_{cnu})/(ker(V_2|\mathfrak{H}^{(1)}_{cnu})^*)]
\oplus\\
&\quad \oplus [(V^*_1|\mathfrak{H}^{(1)}_u)/\ker(V^*_2|\mathfrak{H}^{(1)}_u)^*].
\end{align*}
In the direct sum it is clear that the first operator is a $C_{0\cdot}$-contraction, while the second is unitary since $\ker(V_2|\mathfrak{H}^{(1)}_u)^*$ reduces $V^*_1|\mathfrak{H}^{(1)}_u$. Thus
\[
T_{cnu} = [(V^*_1|\mathfrak{H}^{(1)}_{cnu})/\ker(V_2|\mathfrak{H}^{(1)}_{cnu})^*]^*
\]
is a $C_{\cdot 0}$ contraction.

So the converse statement to the above Corollary is also valid.

2)~~If $\{V_1,V_2\}$ is a bi-isometry and $\dim \ker V^*_2<\infty$, then for the c.n.u.\  part of $T(= (V^*_1|\ker V^*_2)^*$ we have $\|T^n_{cnu}\|\to 0$. Thus the Corollary applies.
\end{rem}

We proceed now to a more detailed analysis of the unitary operator
$Z$. Define a contraction $T'$ on $\mathfrak{F}'$ by
$T'f'=P_{\mathfrak{F}'}Z(0\oplus f')$,
$f'\in\mathfrak{F}'$. The Julia--Halmos matrix associated with $T'$
must again be ``part\char`\"{} of the operator $Z$. More precisely,
consider the decompositions
$\mathfrak{D}_{T}\oplus\mathfrak{F}'=\mathfrak{R}\oplus(\mathfrak{F}'\vee
Z^{*}\mathfrak{F}'$,
$\mathfrak{D}_{T^{*}}\oplus\mathfrak{F}'=\mathfrak{R}_{*}\oplus(\mathfrak{F}'
\vee Z\mathfrak{F}')$,
$\mathfrak{D}_{T}=\mathfrak{R}\oplus\mathfrak{D}'_{*}$, and
$\mathfrak{D}_{T^{*}}=\mathfrak{R}_{*}\oplus\mathfrak{D}'$,
so that $\mathfrak{D}'_{*}\oplus\mathfrak{F}'=\mathfrak{F}'\vee
Z^{*}\mathfrak{F}'$
and $\mathfrak{D}'\oplus\mathfrak{F}'=\mathfrak{F}'\vee Z\mathfrak{F}'$.
We have
$Z(\mathfrak{D}'_{*}\oplus\mathfrak{F}')=\mathfrak{D}'\oplus\mathfrak{F}'$
and $Z\mathfrak{R}=\mathfrak{R}_{*}$. As before, there exist unique
unitary operators $X:\mathfrak{D}_{T'}\to\mathfrak{D}'$ and
$X_{*}:\mathfrak{D}_{T^{\prime*}}\to\mathfrak{D}'_{*}$
such that \[
Z(d'_{*}\oplus
f')=[-XT^{\prime*}X_{*}^{*}d'_{*}+XD_{T'}f']\oplus[D_{T^{\prime*}}X_{*}^{*}d'_{*
}+T'f']\]
 for $d'_{*}\oplus f'\in\mathfrak{D}'_{*}\oplus\mathfrak{F}'$. In
other words, $Z$ is uniquely determined by $T',X,X_{*}$, and the
unitary operator $Y=Z|\mathfrak{R}:\mathfrak{R}\to\mathfrak{R}_{*}$.
We summarize again for future use.

\begin{prop}
Consider a nonet
$(\mathfrak{F},\mathfrak{F}',\mathfrak{R},\mathfrak{R}_{*},T,T',X,X_{*},Y)$,
where $\mathfrak{F}$ and $\mathfrak{F}'$ are Hilbert spaces, $T$
and $T'$ are contractions on $\mathfrak{F}$ and $\mathfrak{F}'$,
respectively, $\mathfrak{R}\subset\mathfrak{D}_{T}$,
$\mathfrak{R}_{*}\subset\mathfrak{D}_{T^{*}}$
are subspaces, and
$X:\mathfrak{D}_{T'}\to\mathfrak{D}_{T^{*}}\ominus\mathfrak{R}_{*}$; and
$X_{*}:\mathfrak{D}_{T^{\prime*}}\to\mathfrak{D}_{T}\ominus\mathfrak{R},$
$Y:\mathfrak{R}\to\mathfrak{R}_{*}$ are unitary operators. Associated
with this data there is a model triple $(\mathfrak{E},U,P)$, where
$\mathfrak{E}=\mathfrak{F}\oplus\mathfrak{D}_{T}\oplus\mathfrak{F}'$,
$P\mathfrak{E}=\mathfrak{F}\oplus\{0\}\oplus\{0\}$, and $U=W_{1}W_{2}$,
with
$W_{1}:\mathfrak{F}\oplus\mathfrak{D}_{T^{*}}\oplus\mathfrak{F}'\to\mathfrak{F}
\oplus\mathfrak{D}_{T}\oplus\mathfrak{F}'$
and $W_{2}:\mathfrak{F}\oplus\mathfrak{D}_{T}\oplus\mathfrak{F}'\to\mathfrak{F}
\oplus\mathfrak{D}_{T^{*}}\oplus\mathfrak{F}'$
given by the matrices
\[
W_{1}=\left[\begin{array}{ccc}
T & D_{T^{*}} & 0\\
D_{T} & -T^{*} & 0\\
0 & 0 & I_{\mathfrak{F}'}\end{array}\right],\quad W_{2}=\left[\begin{array}{ccc}
I_{\mathfrak{F}} & 0&0\\
0\\
0 & \multicolumn{2}{c}{Z}\end{array}\right],\]
 where $Z:\mathfrak{D}_{T}\oplus\mathfrak{F}'=\mathfrak{R}\oplus
X_{*}\mathfrak{D}_{T^{\prime*}}\oplus\mathfrak{F}'\to\mathfrak{D}_{T^{*}}\oplus
\mathfrak{F}'=\mathfrak{R}_{*}\oplus X\mathfrak{D}_{T'}\oplus\mathfrak{F}'$
is the unitary given by\[
Z=\left[\begin{array}{ccc}
Y & 0 & 0\\
0 & -XT^{\prime*}X_{*}^{*} & XD_{T'}\\
0 & D_{T^{\prime*}}X_{*}^{*} & T'\end{array}\right].\]

\begin{enumerate}
\item Every model triple can be obtained, up to unitary equivalence, in
the manner described above.
\item Two nonets
$(\mathfrak{F}_{j},\mathfrak{F}_{j}',\mathfrak{R}_{j},\mathfrak{R}_{*j},T_{j},T_
{j}',X_{j},X_{*j},Y_{j}),$
$j=1,2,$ determine equivalent model triples if and only if there
exist unitary operators $A:\mathfrak{F}_{1}\to\mathfrak{F}_{2}$ and
$A':\mathfrak{F}_{1}'\to\mathfrak{F}'_{2}$ satisfying $AT_{1}=T_{2}A$,
$A'T'_{1}=T_{2}'A'$, $A\mathfrak{R}_{1}=\mathfrak{R}_{2}$,
$A'\mathfrak{R}_{*1}=\mathfrak{R}_{*2}$,
$(A'|\mathfrak{R}_{*1})Y_{1}=Y_{2}(A|\mathfrak{R}_{1})$,
$AX_{1}=X_{2}A'|\mathfrak{D}_{T_{1}}$,
and $AX_{*1}=X_{*2}A'|\mathfrak{D}_{T_{1}^{\prime*}}$.
\item Two contractions $T$ and $T'$ can be included in one of the nonets
described above if and only if
$\dim\mathfrak{D}_{T'}\le\dim\mathfrak{D}_{T^{*}}$,
$\dim\mathfrak{D}_{T^{\prime*}}\le\dim\mathfrak{D}_{T}$, and \[
\dim\mathfrak{D}_{T}+\dim\mathfrak{D}_{T'}=\dim\mathfrak{D}_{T^{*}}+\dim
\mathfrak{D}_{T^{\prime*}}.\]

\end{enumerate}
\end{prop}
Part (3) of the above statement simply gives conditions for the existence
of subspaces $\mathfrak{R}$ and $\mathfrak{R}_{*}$ such that the
various unitaries in a nonet can be constructed.

\section{Irreducible Multi-Isometries}

In the classification of multi-isometries, it seems natural to consider
first the irreducible ones, that is, those which do not have a common
reducing subspace. Assume that $(V_{1},V_{2},\dots,V_{n})$ is an
irreducible $n$-isometry, and $V_{j}$ is a unitary operator for
some $j$. The spectral subspaces of $V_{j}$ are hyperinvariant for
$V_{j}$, hence reducing for the $n$-isometry. We conclude that $V_{j}$
has no nontrivial spectral subspaces, so that $V_{j}=\lambda_{j}I$
for some scalar $\lambda_{j}$. Clearly our $n$-isometry will be
as easy to study as the $(n-1)$-isometry obtained by deleting $V_{j}$.
An $n$-isometry will be said to be proper if none of the component
isometries is a constant multiple of the identity operator.

\begin{lem}
Any irreducible proper $n$-isometry $(V_{1},V_{2},\dots,V_{n})$
is cnu, and $V_{1}V_{2}\cdots V_{n}$ is a unilateral shift of multiplicity
at least $n$.
\end{lem}
\begin{proof}
The unitary part is a reducing space for the $n$-isometry, so it
must be trivial. As noted above, none of the $V_{j}$ is unitary,
and therefore all the inclusions \[
\mathfrak{H}\supset V_{1}\mathfrak{H}\supset
V_{1}V_{2}\mathfrak{H}\supset\cdots\supset V_{1}V_{2}\cdots V_{n}\mathfrak{H}\]
 are strict, thus proving the last assertion of the lemma.
\end{proof}
We will see that all the $V_{j}$ must in fact be pure isometries if
$(V_{1},V_{2},\dots,V_{n})$
is an irreducible proper isometry such that $V_{1}V_{2}\cdots V_{n}$
has finite multiplicity. We need a preliminary result, which follows from Corollary 3.6.

\begin{lem}
Let $V_{1}$ and $V_{2}$ be commuting isometries such that $\dim\ker
V_{2}^{*}<\infty$.
Consider the von Neumann--Wold decomposition $V_{1}=S\oplus U$ on
$\mathfrak{H}\oplus\mathfrak{H}'$
such that $S$ is pure and $U$ is unitary. Then $\mathfrak{H}$ is
a reducing subspace for $V_{2}$.
\end{lem}

This lemma enables us to prove the following

\begin{cor}
Let $F$ be an irreducible family of commuting isometries such that
$\dim\ker V^{*}<\infty$ for every $V\in F$. Then each $V\in F$
is either pure, or a constant multiple of the identity.
\end{cor}
\begin{proof}
Assume that a $V\in F$ is neither pure, nor unitary. Then the preceding
lemma provides a reducing subspace for $F$. If $V$ is unitary, but
not a scalar multiple of the identity, then any nontrivial spectral
space for $V$ reduces $F$.
\end{proof}
\begin{cor}
Let $(V_{1},V_{2},\dots,V_{n})$ be an irreducible proper $n$-isometry
such that $V_{1}V_{2}\cdots V_{n}$ is a pure isometry of multiplicity
$\le2n-1$. Then each $V_{i}$ is pure, and one of them has multiplicity
one. When $V_{1}$ has multiplicity one, we must have $V_{j}=\varphi_{j}(V_{1})$
where $\varphi_{j}$ is a finite Blaschke product for $j=2,3,\dots,n$,
and the sum of the multiplicities of the $\varphi_{j}$ must be at most $2n-2$.
If $V_{1}V_{2}\cdots V_{n}$ has multiplicity $n$, then each $V_{i}$
must have multiplicity one.
\end{cor}
\begin{proof}
The preceding result shows that the  $V_{i}$ are pure, so that the
inclusions \[
\mathfrak{H}\supset V_{1}\mathfrak{H}\subset
V_{1}V_{2}\mathfrak{H}\supset\cdots\supset V_{1}V_{2}\cdots V_{n}\mathfrak{H}\]
 are strict, and \[
\dim[(V_{1}V_{2}\cdots V_{i-1}\mathfrak{H})\ominus(V_{1}V_{2}\cdots
V_{i}\mathfrak{H})]=\dim[\mathfrak{H}\ominus V_{i}\mathfrak{H}],\quad
i=1,2,\dots,n.\]
 The conclusion follows because the sum of these $n$ positive integers
is $\le2n-1$.
\end{proof}
\begin{cor}
Let $(V_{1},V_{2},\cdots,V_{n})$ be an irreducible proper $n$-isometry,
and assume that $V_{1}V_{2}\cdots V_{n}$ has multiplicity $n$. Then
there exist M\"obius transformations $\varphi_{j}$
such that $V_{j}=\varphi_{j}(V_{1})$ for $j=2,3,\dots,n.$
\end{cor}
\begin{proof}
The fact that  $V_{j}=\varphi_{j}(V_{1})$ for some $\varphi_{j}\in H^{\infty}$
follows because $V_{j}$ commutes with $V_{1}$, and $V_{1}$ is a
pure isometry of multiplicity one. Moreover $\varphi(V_{1})$ is an
isometry if and only if $\varphi$ is an inner function, and the multiplicity
of $\varphi(V_{1})$ is the multiplicity of $\varphi$.
\end{proof}
The last result allows us to classify completely all irreducible proper
$n$-isometries $(V_{1},V_{2},\cdots,V_{n})$ for which $V_{1}V_{2}\cdots V_{n}$
is a shift of multiplicity $n$. Indeed, up to unitary equivalence,
we may assume that $V_{1}$ is the shift $S$ of multiplicity one
on $H^{2}$, so that our $n$-isometry is
$(S,\varphi_{2}(S),\dots,\varphi_{n}(S))$
for some M\"obius transforms $\varphi_{2},\dots,\varphi_{n}$.

\begin{prop}
Let $(\varphi_{j})_{j=2}^{n}$ and $(\psi_{j})_{j=2}^{n}$ be two
families of inner functions. The $n$-isometries
$(S,\varphi_{2}(S),\dots,\varphi_{n}(S))$
and $(S,\psi_{2}(S),\dots,\psi_{n}(S))$ are unitarily equivalent
if and only if $\varphi_{j}=\psi_{j}$ for $j=2,\dots,n$.
\end{prop}
\begin{proof}
Let $U$ be a unitary operator on $H^{2}$ satisfying $US=SU$ and
$U\varphi_{j}(S)=\psi_{j}(S)U$. Then $U$ must in fact be a scalar
multiple of the identity, so that $\varphi_{j}(S)=\psi_{j}(S)$ and
therefore $\varphi_{j}=\psi_{j}$.
\end{proof}
This result classifies $n$-isometries whenever $V_{1}$ is a shift
of multiplicity one.
\medskip

As mentioned in the introduction, the unitary invariants of irreducible
proper $n$-isometries, such that $V_{1}V_{2}\cdots V_{n}$ has multiplicity
$n$, can be described explicitly for $n=2$ and $n=3$. When $n=2$
we must simply describe (up to unitary equivalence) all pairs $(U,P)$,
where $U$ is unitary and $P$ is a projection of rank one on a space
$\mathfrak{E}$ of dimension 2. Using unitary equivalence, we may
assume that $\mathfrak{E}=\mathbb{C}^{2}$, and \[
P=\left[\begin{array}{cc}
0 & 0\\
0 & 1\end{array}\right].\]
 The possible unitary operators $U$ are given by \[
U(c,\theta)=\left[\begin{array}{cc}
c & d\theta\\
d &
\bar{c}\theta\end{array}\right]\:\text{with}\:|\theta|=1,|c|<1,\:\text{and}\:
d=(1-|c|^{2})^{1/2}.\]
 Moreover, if $W$ is a unitary operator on $\mathbb{C}^{2}$ such
that $WP=PW$ and $WU(c,\theta)=W(c'\theta')U$, the reader will verify
with no difficulty that necessarily $c=c'$ and $\theta=\theta'$. Thus
$(c,\theta)$ is a complete set of unitary invariants for 2-isometries
$(V_1,V_2)$ for which the multiplicity of $V_1V_2$ is two.

When $n=3$, the unitary invariants consist of commuting unitaries
$U_{1},U_{2}$ and projections $P_{1},P_{2}$ of rank one on a space
$\mathfrak{E}$ of dimension 3 satisfying
$P_{1}+U_{1}^{*}P_{1}U_{1}=P_{2}+U_{2}^{*}P_{2}\le I_{\mathfrak{E}}$.
As before, we may assume that $\mathfrak{E}=\mathbf{C}^{3}$, \[
P_{1}=\left[\begin{array}{ccc}
0 & 0 & 0\\
0 & 0 & 0\\
0 & 0 & 1\end{array}\right]\:\text{and}\;
U_{1}^{*}P_{2}U_{1}=\left[\begin{array}{ccc}
0 & 0 & 0\\
0 & 1 & 0\\
0 & 0 & 0\end{array}\right].\]
 Denoting by $e_{j}$ the standard basis vectors in $\mathbb{C}^{3}$,
there will exist complex numbers $\alpha,\beta,\gamma,\delta,\varepsilon,\eta$
such that \[
U_{2}^{*}e_{3}=\overline{\alpha}e_{3}+\beta e_{2},U_{1}^{*}U_{2}^{*}e_{3}=\gamma
e_{3}+\delta e_{1},U_{1}e_{2}=\varepsilon e_{3}+\eta e_{2}.\]
Setting $N=(U_{1}-\eta)(U_{2}^{*}-\overline{\alpha})-\beta\varepsilon,$
we note that $N$ is normal and $Ne_{2}=Ne_{3}=0.$ If $N\ne0$, it
follows that the span of $e_{2}$ and $e_{3}$ reduces $U_{1},U_{2},P_{1},$
and $P_{2},$ so that the corresponding 3-isometry is reducible as
well. Thus in the irreducible case we must have  $N=0$. An easy (but
tedious) calculation shows that irreducibility also implies
$|\alpha|<1,|\eta|<1,$
and $\beta\varepsilon\ne0$. The equation $N=0$ can then be solved
for $U_{1}$, yielding\[
U_{1}=\theta(U_{2}-\alpha)(I-\overline{\alpha}U_{2})^{-1},\]
where the number $\theta=\varepsilon/\overline{\beta}$ has absolute
value one. Finally, applying a unitary equivalence with a diagonal
operator on $\mathbb{C}^{3}$, we can assume that \[
U_{3}=\left[\begin{array}{ccc}
\alpha_{1} & \theta_{1}\overline{\alpha}d_{1} & -\theta_{1}dd_{1}\\
d_{1} & -\theta_{1}\overline{\alpha\alpha_{1}} &
\theta_{1}\overline{\alpha_{1}}d\\
0 & d & \alpha\end{array}\right],\]
where $d=(1-|\alpha|^{2})^{1/2},d_{1}=(1-|\alpha_{1}|^{2})^{1/2},$
and $|\theta_{1}|=1.$ We must also have $|\alpha_{1}|<1$ on account
of irreducibility. Given numbers $\alpha,\alpha_{1},\theta,\theta_{1}$
such that $|\alpha|<1,$ $|\alpha_{1}|<1,$ and $|\theta|=|\theta_{1}|=1,$
let us denote by $U_{1}(\alpha,\alpha_{1},\theta,\theta_{1})$ and
$U_{2}(\alpha,\alpha_{1},\theta,\theta_{1})$ the unitary operators
given by the above formulas. In this way the quadruple
$(\alpha,\alpha_1,\theta,\theta_1)$ determines a 3-isometry $(V_1,V_2,V_3)$ such
that the multiplicity of $V_1V_2V_3$ is three. Again, the reader will be able to
verify
that a unitary $W$ satisfying $WP_{j}=P_{j}W$ and
$WU_{1}(\alpha,\alpha_{1},\theta,\theta_{1})=U_{1}(\alpha',\alpha'_{1},\theta',
\theta'_{1})W$
exists only in case $\alpha=\alpha'$ , $\alpha_{1}=\alpha_{1}^{'},$
$\theta=\theta'$, and $\theta_{1}=\theta_{1}'$.

It is interesting to note that, in case $V_{1}V_{2}\cdots V_{n}$
is a shift of multiplicity at least $2n$ , the isometries $V_{j}$
need not all belong to the algebra generated by some shift of multiplicity
one. A simple example is obtained for $n=2$ with $V_{1}=S\oplus S$,
and \[
V_{2}=\left[\begin{array}{cc}
0 & S^{2}\\
I & 0\end{array}\right],\]
where $S$ denotes the standard shift of multiplicity one on $H^{2}$.
For this example we have $V_{1}^{2}=V_{2}^{2},$ so that
$(V_{1}-V_{2})(V_{1}+V_{2})=0,$
while $V_{1}-V_{2}\ne0\ne V_{1}+V_{2}$. The commutant of a shift
of multiplicity one is isomorphic to $H^{\infty}$, and this algebra
has no zero divisors; thus $V_{1}$ and $V_{2}$ cannot belong to
the commutant of the same shift of multiplicity one.

\end{document}